\documentclass[sn-mathphys]{sn-jnl}

\usepackage{xspace}
\usepackage{makecell} 
\usepackage{amsmath,amsfonts,amssymb,amsthm}
\usepackage{mathtools}
\usepackage{tabularx,booktabs}
\usepackage{multirow}
\usepackage[justification=centering]{caption}

\jyear{2021}%

\theoremstyle{thmstyleone}%
\theoremstyle{thmstyletwo}%

\theoremstyle{thmstylethree}%

\def\KS{KS\xspace}
\def\KSF{KS'\xspace}
\def\KSFVG{KS(V, G)\xspace}
\def\KSFV{KS(V)\xspace}

\def\CPXA{CPX-A\xspace}
\def\CPXB{CPX-B\xspace}
\def\CPXC{CPX-C\xspace}
\def\VMCP{VMCP\xspace}
\def\VMCPs{VMCPs\xspace}

\def\NVM{\#VM\xspace}

\def\TT{T\xspace}

\def\GP{Gap\xspace\%\xspace}
\def\WGP{Worst Gap\xspace\%\xspace}

\def\VM{VM\xspace}
\def\VMs{VMs\xspace}

\def\CASSIGN{c^{\rm alloc}_{i,j}}
\def\CRUN{c^{\rm run}_{j}}
\def\CMIG{c^{\rm mig}_{i,j}}
\def\CNEW{c^{\rm new}_{i,j}}
\def\UIR{u_{i,r}}
\def\SKR{s_{j,r}}
\def\DI{d_{i}}
\def\DINEW{d^{\rm new}_{i}}
\def\NIK{n_{i,j}}

\def\ISET{\mathcal{I}}
\def\KSET{\mathcal{J}}
\def\RSET{\mathcal{R}}
\def\XIK{x_{i,j}}
\def\YK{y_{j}}
\def\ZIK{z_{i,j}}
\def\XIKNEW{x^{\rm new}_{i,j}}

\newcommand{\MIP}{{MILP}\xspace}
\newcommand{\LP}{{LP}\xspace}
\newcommand{\solver}[1]{\textsc{#1}\xspace}
\newcommand{\cplex}{\solver{CPLEX}}

\begin{document}

\title[A Kernel Search Algorithm for Virtual Machine Consolidation Problem]{A Kernel Search Algorithm for Virtual Machine Consolidation Problem}


\author[1]{\fnm{Jiang-Yao} \sur{Luo}}\email{luoshui3000@bupt.edu.cn}

\author*[1]{\fnm{Jian-Hua} \sur{Yuan}}\email{jianhuayuan@bupt.edu.cn}


\affil*[1]{\orgdiv{School of Science}, \orgname{Beijing University of Posts and Telecommunications}, 
	\city{Beijing}, \postcode{100876}, 
		\country{China}}



\abstract{
	Virtual machine consolidation describes the process of reallocation of virtual machines (\VMs) on a set of target servers.
	It can be formulated as a mixed integer linear programming problem which is proven to be an NP-hard problem. 
	In this paper, we propose a kernel search (\KS) heuristic algorithm based on hard variable fixing to quickly obtain a high-quality solution for large-scale virtual machine consolidation problems (\VMCPs). 
	Since variable fixing strategies in existing \KS works may make \VMCP infeasible, our proposed \KS algorithm employs a more efficient strategy to choose a set of fixed variables according to the corresponding reduced cost.
	Numerical results on \VMCP instances demonstrate that our proposed \KS algorithm significantly outperforms the state-of-the-art mixed integer linear programming solver in terms of CPU time, and our proposed strategy of variable fixing significantly improves the efficiency of the \KS algorithm as well as the degradation of solution quality can be negligible.
}




\maketitle

\section{Introduction}
\label{sect:introduction}
Nowadays, due to the extensibility and flexibility of cloud computing, 
more and more internet services are provided by it.
Indeed, internet services can be instantiated inside virtual machines (\VMs) and flexibly allocated to any servers in the cloud data center. 
However, as the virtual machines dynamically change, resulting in the improper allocation of \VMs and imbalanced load distribution, 
the high energy consumption and low efficiency of cloud data centers are also becoming more and more serious. 
In fact, for large data centers, 15 to 20 percent resource utilization is common \cite{vogels2008beyond}, 
and even if an activated server is kept idle, it consumes up to 66 percent of the peak power \cite{chen2008energy}.
To improve efficiency and reduce energy consumption, cloud service providers employ \VM migration technology to dynamically reallocate \VMs 
by migrating old \VMs among servers and mapping new \VMs into servers in data centers.
The above problem is called the \emph{virtual machine consolidation problem} in the literature.

In general, virtual machine consolidation can be formulated as a mixed integer linear programming (\MIP) problem 
which determines the activated servers and the reallocation of \VMs in such a way that 
the sum of server activation, \VM allocation, and \VM migration costs is minimized 
subject to resource constraints of the servers and other practical constraints.
The \VMCP is strongly NP-hard \cite{speitkamp2010mathematical}, 
so there is no polynomial time algorithm to solve the \VMCP to optimality unless P=NP.
Therefore, existing works mainly focus on  
heuristic \cite{speitkamp2010mathematical,mazumdar2017power,beloglazov2012optimal,verma2008pmapper,goudarzi2012sla} 
and metaheuristic \cite{sharma2016multi,li2020energy,wu2016energy,jiang2017dataabc} 
algorithms for solving the \VMCP.
Among heuristic algorithms, greedy heuristics are the popular algorithm used to solve \VMCPs \cite{wei2020exact,varasteh2015server}, which are established 
based on the original first fit, best fit, first fit decreasing, and best fit decreasing mechanisms \cite{beloglazov2012optimal,verma2008pmapper}.
Greedy heuristics are used as a baseline for comparison in many works because they can provide fast solutions, 
but they are more problem-dependent than other algorithms.
More specifically, their design always depends on the specific problem structure, thus, 
is not easily extendable to other similar problems.
References \cite{speitkamp2010mathematical} and \cite{goudarzi2012sla} proposed some more problem-independent heuristics than greedy algorithms. 
However, there is still a non-trivial gap, from 6\% to 49\% on average, between the solutions found by the heuristics with the optimal solution \cite{speitkamp2010mathematical,goudarzi2012sla,mazumdar2017power}.
Although metaheuristic algorithms \cite{sharma2016multi, wu2016energy,li2020energy,jiang2017dataabc} are problem-independent and can provide better solutions than greedy algorithms,
they require more iteration time.
Fortunately, the kernel search (\KS) algorithm is problem-independent and can quickly obtain a high-quality solution.
The standard \KS, first proposed by reference \cite{angelelli2010kernel}, has been applied to the solution of 
the multi-dimensional knapsack problem that is strongly NP-hard.  
Subsequently, the \KS based algorithms have been successfully applied to solve 
the portfolio optimization problems \cite{angelelli2012kernel},
facility location problems \cite{guastaroba2012facility, guastaroba2014heuristic, filippi2021kernel}, 
index tracking problems \cite{guastaroba2012index, filippi2016heuristic}, 
and general mixed integer programming problems \cite{guastaroba2017adaptive}.
Moreover, the \KS heuristic requires little implementation effort 
since the most cumbersome part of the search is finished by a state-of-the-art \MIP solver \cite{angelelli2012kernel}.
Unfortunately, to the best of our knowledge, there is no research on using \KS for solving \VMCP.
Furthermore, our preliminary experiments show that the standard \KS for solving \VMCP is time-consuming.
Therefore, a variant of kernel search based on hard variables fixing is presented to address \VMCP.
A crucial issue in variable fixing is related to the strategy of variable fixing, 
which significantly determines the efficiency of the \KS algorithm.
Reference \cite{guastaroba2014heuristic} proposed a variant of \KS based on hard variable fixing 
that only fixes some binary variables to their values in \LP solution.
However, our preliminary experiments indicate that this strategy of variable fixing 
results in some \VMCP instances being infeasible.
For this reason, our proposed variant adopts a more sophisticated but efficient strategy of variable fixing 
that chooses the fixed binary variables according to corresponding reduced costs
to avoid infeasibility occurring.
As integer variables dominate the \MIP formulation of \VMCP,
so we apply a similar strategy of variable fixing for them.
For more details introduction to hard variable fixing, we refer readers to references \cite{danna2005exploring,lazic2010variable}.
The main contribution of this paper is that we propose a \KS heuristic algorithm based on hard variable fixing, 
and apply it to quickly obtain a high-quality solution for the large-scale \VMCPs.
We provide a new strategy of variable fixing to enhance the efficiency of the \KS algorithm and avoid \VMCP infeasibility occurring.
In addition to fixing binary variables, we also fixed the integer variables which dominate the \VMCP.
Extensive computational results show that our proposed \KS algorithm significantly outperforms three settings of the standard \MIP solver that emphasizes the heuristics, and our proposed strategy of variable fixing significantly improves the efficiency of the \KS algorithm as well as the degradation of solution quality can be negligible.

\indent The paper is organized as follows. 
Section \ref{sect:problem_formualtion} presents the \MIP formulation of the \VMCP. 
Section \ref{sect:solution_methodology} describes the proposed kernel search algorithm to solve \VMCPs. 
Section \ref{sect:numerical_results} shows computational results. 
Finally, Section \ref{sect:conclusion} draws some concluding remarks. 

\section{Virtual machine consolidation problem}
\label{sect:problem_formualtion}
Virtual machine consolidation describes the process of combining several different virtual machines and assigning them to a set of target servers.
It can be used to optimize the allocation of \VMs and servers for minimizing the allocation costs of \VMs, the activation cost of servers, and the migration costs of \VMs.
In this section, we follow \cite{luo2022cutandsolve} to present a compact mixed integer linear programming formulation for the \VMCP.

\begin{figure}[!htbp]
	\centering
	\includegraphics[width=0.8\linewidth]{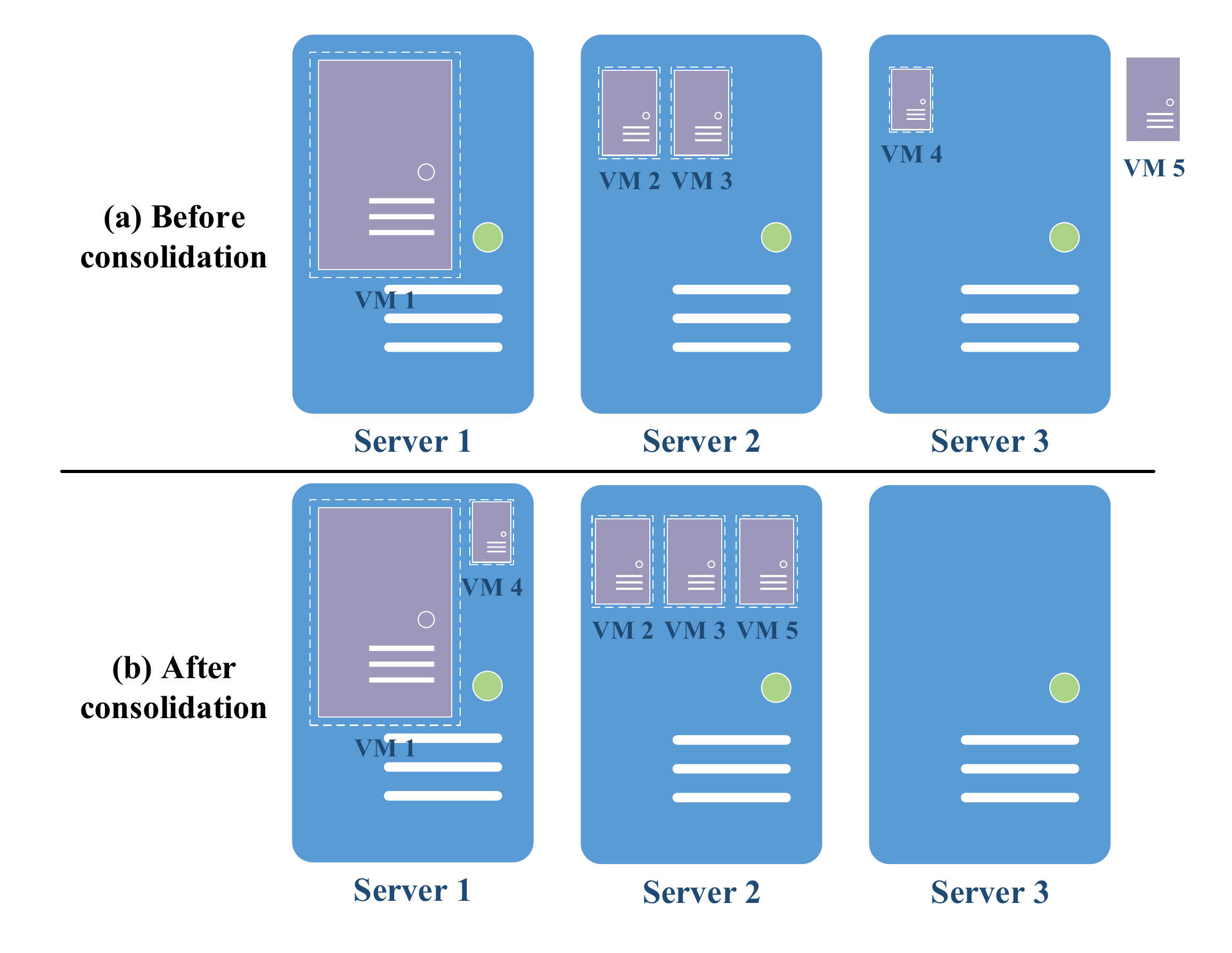}
	\caption{An example of the virtual machine consolidation problem} 
	\label{fig:vmcp}
\end{figure}
Fig. \ref{fig:vmcp} depicts an example of the virtual machine before and after consolidation.
Let $\KSET$, $\ISET$, and $\RSET$ denote the set of the servers,
the set of types of \VMs that needs to be allocated to the servers,
and the set of resources of the servers, respectively.
As shown in Fig. \ref{fig:vmcp}, we are concerned with a set of servers associated with 
a certain capacity $\SKR$ for each $j \in \mathcal{J}, r \in \mathcal{R}$, and a set of 
old and new virtual machines associated with the demand $\UIR$ for each $i \in \mathcal{I}, r \in \mathcal{R}$.
Before \VM consolidation, there are $\NIK$ old \VMs of type $i$ (e.g., \VMs 1-4 in Fig. \ref{fig:vmcp}) that are currently allocated at server $j$
and $\sum_{i \in \ISET}\DINEW$ new incoming \VMs (e.g., \VM 5 in Fig. \ref{fig:vmcp}) that are needed to be allocated at servers.
For notation purpose, we denote $\DI = \sum_{j \in \KSET} \NIK$ for all $i \in \ISET$.

After \VM consolidation, we introduce integer variable $\XIK$ to indicate  the number of old \VMs of type $i$ allocated to server $j$, 
binary variable $\YK$ to represent whether or not server $j$ is activated,
integer variable $\ZIK$ to indicate the number of old \VMs of type $i$ migrated to server $j$, 
and integer variable $\XIKNEW$ to denote the number of new incoming \VM of type $i$ allocated to server $j$.
As discussed, the \VMCP can be formulated as follows: 
\begin{subequations}
\label{eq:mathforms}
\begin{align} 
\!\!\!\! \mathop{\mathrm{min}} \ 
&\sum_{i \in \ISET}\sum_{j \in \KSET} \CASSIGN \NIK
 +\sum_{j \in \KSET} \CRUN \YK  + \sum\limits_{i \in \ISET}\sum\limits_{j \in \KSET} \CMIG \ZIK 
+\sum\limits_{i \in \ISET}\sum\limits_{j \in \KSET}\CNEW \XIKNEW  \label{eq:obj}\\
\!\!\!\! \text{s.t.} ~~& \sum\limits_{i \in \ISET} \UIR \XIK + \sum\limits_{i \in \ISET}\UIR\XIKNEW \leq \SKR \YK,~  
\forall~r \in \RSET,~\forall~j \in \KSET,  \label{eq:resource}   \\
&  \XIK - \NIK \leq \ZIK , ~ \forall~i \in \ISET,~\forall~j \in \KSET,  \label{eq:migration} \\
& 	\sum\limits_{k \in \KSET} \XIK = \DI, ~ \forall~i \in \ISET,   \label{eq:allocation_general}  \\
& \sum\limits_{k \in \KSET} \XIKNEW = \DINEW,~\forall~i \in \ISET, \label{eq:new_incoming_vms} \\
& \XIK, \ZIK, \XIKNEW \in \mathbb{Z_+},~\XIK \leq v_{i,j},~\forall~i \in \ISET,~\forall~j \in \KSET, \label{eq:varX}\\
&  \YK \in \{0,1\},~\forall~j \in \KSET. \label{eq:varY}
\end{align}
\end{subequations}
The objective function \eqref{eq:obj} to be minimized is the cost of allocating all old \VMs to servers, 
the sum of the  activation cost of servers, 
the cost of migrating \VMs among servers,
and the cost of assigning all new incoming \VMs to servers. 
Here coefficients $\CRUN$, $\CASSIGN$, $\CMIG$, and $\CNEW$ are greater or equal to zero, which represent the activation cost of server $j$, 
the cost of allocating a old \VM of type $i$ to server $j$, the cost of migrating a \VM of type $i$ to server $j$, 
and the cost of assigning a new incoming \VM of type $i$ to server $j$, respectively. 

Constraint \eqref{eq:resource} ensures that any type of resource capacity for each server is enough for
the aggregate workload of \VMs allocated 
Constraint \eqref{eq:migration} makes sure that if the number of \VMs of type $i$ allocated to server $j$ 
after the consolidation, $\XIK$, is larger than 
that before the consolidation, $\NIK$, then the number of \VMs of type $i$ migrated to server $j$, $\ZIK$,  
must be equal to $\XIK-\NIK$ due to $\CMIG \geq 0$; 
otherwise, it is equal to zero due to $\ZIK \in \mathbb{Z_+}$.
Constraints \eqref{eq:allocation_general} and \eqref{eq:new_incoming_vms} state that old and new \VMs of each type have to be assigned to servers, respectively. 
Finally, constraints \eqref{eq:varX} and \eqref{eq:varY} restrict $\XIK$, $\ZIK$, $\XIKNEW$, and $\YK$ to be integer/binary variables and trivial upper bounds $\left\{\XIK\right\}$ for variables $\left\{v_{i,j}\right\}$ where
\begin{equation*}
	v_{i,j} = \min\left\{\sum_{j \in \KSET}n_{i,j}, ~\min_{r \in \RSET} \left \{ \left\lfloor  \frac{s_{j,r}}{u_{i,r}} \right\rfloor \right  \}\right\}, ~\forall~i \in \ISET, ~\forall~j \in \KSET.
\end{equation*} 
In the next section, we shall develop an efficient customized algorithm based on the kernel search heuristic for solving the large-scale \VMCPs.
\section{The kernel search algorithm}
\label{sect:solution_methodology} 
In this section, we shall first provide an intuitive and general description of the standard kernel search algorithm 
and then develop a new kernel search algorithm which is designed 
to quickly obtain a high-quality feasible solution for the large-scale \VMCPs. 
Most of the notation and definitions introduced in this section will be used throughout the remainder of the paper. 

\subsection{The standard kernel search algorithm}
\label{sect:standard_kernel_search}
Our description of the standard \KS algorithmic framework mainly refers to references
 \cite{angelelli2010kernel,angelelli2012kernel,guastaroba2014heuristic, guastaroba2017adaptive}.
To simplicity the notation, let $\mathcal{V}$ and $\mathcal{G}$ denote the set of binary variables and 
the set of general integer variables in \VMCP \eqref{eq:mathforms}, respectively. 
We refer to the \MIP problem including all variables in $\mathcal{V} \cup \mathcal{G}$ as the original problem, 
and call the restricted problem where the binary variables in $\mathcal{V} \backslash \mathcal{U}$ ($\mathcal{U} \subsetneqq \mathcal{V}$) are fixed to zero by \MIP($\mathcal{U}$).

\begin{figure}[!htbp]
	\centering
	\includegraphics[width=\linewidth]{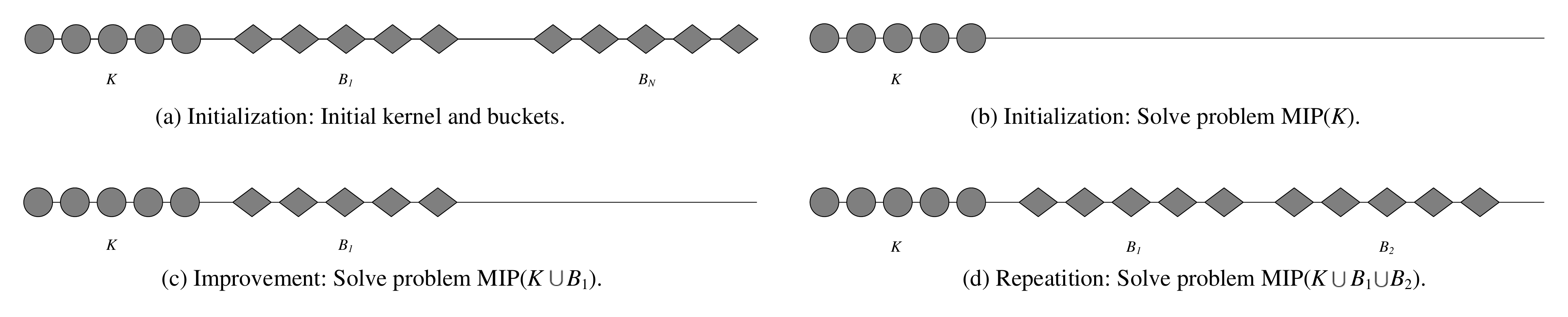}
	\caption{An illustrative example of the operation of the standard \KS algorithm.}
	\label{fig:ks_framework}
\end{figure}
\KS is essentially a heuristic framework with a general and flexible structure applicable to any \MIP problem with binary variables.
Fig. \ref{fig:ks_framework} describes an iteration in the standard \KS algorithmic framework.
More specifically, using the information provided by the optimal solution of the (root node) linear programming (\LP) relaxation of the original problem, 
the standard \KS framework generates possibly different orders for binary variables  
according to the design criterion in the initialization phase.
According to common design criterion, the more left the binary variable in Fig. \ref{fig:ks_framework}(a) is ranked, 
the more likely it is that the value of 1 is taken in the optimal solution of the original problem.
We select the first $\vert \mathcal{K} \rvert$ promising binary variables (see the circles in Fig. \ref{fig:ks_framework}(a))  
from the left to construct the initial kernel $\mathcal{K} \subset \mathcal{V}$.
The remaining binary variables in $\mathcal{V} \backslash \mathcal{K}$ (see the rhombus in Fig. \ref{fig:ks_framework}(a)) 
are partitioned into $N$ subsets (called buckets) denoted as $\mathcal{B}_{i},~i =1,\cdots,N$.
Then we first solve the restricted problem by only considering variables in $\mathcal{K} \cup \mathcal{G}$ 
(see problem \MIP($\mathcal{U}$) where $\mathcal{U} \coloneqq \mathcal{K}$ in Fig. \ref{fig:ks_framework}(b)), 
denoted its solution by $(v_{\rm min},g_{\rm min})$ and the corresponding objective value by ${\rm UB}_{\rm min}$.
\begin{itemize}
	\item[(i)] If set $\mathcal{K}$ is too small, the optimal solution quality of \MIP($\mathcal{U}$) where $\mathcal{U} \coloneqq \mathcal{K}$ is poor or even infeasible (solution quality);
	\item[(ii)] If set $\mathcal{K}$ is too large, we cannot find the optimal solution of \MIP($\mathcal{U}$) where $\mathcal{U} \coloneqq \mathcal{K}$ within a reasonable time (solution efficiency).
\end{itemize}
Neither (i) nor (ii) is our ideal situations. 
If total computational time after solving the first restricted problem is still within the predefined time limit $T_{\rm max}$, 
we solve other restricted problems (see Figs. \ref{fig:ks_framework}(c) and \ref{fig:ks_framework}(d)) 
in the sequence considering the previous $\mathcal{U}$ plus the binary variables belonging to bucket $\mathcal{B}_{i}$.
Two additional constraints 
\begin{equation}
\eqref{eq:obj} \leq {\rm UB}_{\rm min}, \label{eq:cutoff}
\end{equation}
and
\begin{equation}
\sum\limits_{j \in B_{i}}v_{j} \geq 1. \label{eq:piercing_cut}
\end{equation}
are introduced before solving restricted problem \MIP($\mathcal{U}$)
to reduce the computational time required by the solver to find the optimal solution.
The above procedure is repeated until 
the number of buckets already analyzed in $\mathcal{U}$ reaches the limit $\bar{N} \leq N$.
The details of standard \KS are summarized in Algorithm \ref{alg:kernel_search}.
\begin{algorithm}[!htbp]
	\caption{The standard \KS algorithm} 
	\label{alg:kernel_search}
	\begin{algorithmic}[1]
		\Require \VMCP, $\bar{N}$, and $T_{\rm max}$.
		\Ensure $\text{UB}_{\min}$.
		\State Initialize $\text{UB}_{\min} \coloneqq +\infty$ and $i \coloneqq  1$;
		\State Solve the (root) \LP relaxation of the original \VMCP; 
		\State Sort the binary variables according to a predefined sorting criterion; 
		\State Construct an initial kernel $\mathcal{K}$ by selecting the first $\vert \mathcal{K} \rvert$ binary variables;
		\State Consider the binary variables belonging to $\mathcal{V} \backslash \mathcal{K}$ and construct a sequence $\left\{\mathcal{B}_{i}\right\}_{i =1,\cdots,N}$ of buckets;
		\State Set maximum time limit $T_{i} \coloneqq T_{\rm max} / (N+1)$ of each bucket; 
		\State Solve problem \MIP($\mathcal{U}$) where $\mathcal{U} \coloneqq \mathcal{K}$;
		\While{$i \leq \bar{N}$}
			\State Add the two following constraints in \MIP($\mathcal{U} \cup \mathcal{B}_{i}$):
			\State (i) \eqref{eq:obj} $\leq {\rm UB}_{\rm min}$;
			\State (ii) $\sum\limits_{j \in \mathcal{B}_{i}} v_{j} \geq 1$;
			\State Solve problem \MIP($\mathcal{U}$) where $\mathcal{U} \coloneqq \mathcal{U} \cup \mathcal{B}_{i}$ within time limit $T_{i}$ and denote its objective value by ${\rm UB}_{i}$; 
			\If{${\rm UB}_{i} < {\rm UB}_{\rm min}$}
			\State Update ${\rm UB}_{\rm min} \coloneqq {\rm UB}_{i}$;
			\EndIf
			\State Set $i \coloneqq i + 1$; 
		\EndWhile
	\end{algorithmic}
\end{algorithm}
\subsection{The proposed kernel search algorithm}
\label{sect:kernel_search_with_fixed_variable}
The kernel search algorithm has been demonstrated to obtain high-quality solutions for various \MIP problems with binary variables
 \cite{angelelli2012kernel,guastaroba2012index,guastaroba2012facility,guastaroba2014heuristic,guastaroba2017adaptive,filippi2021kernel}.
However, our preliminary experiments showed that due to the solution space of the restriction problem for \KS is still too large, 
the restricted problem is time-consuming, leading to very poor performance of \KS for \VMCPs.
Fortunately, the efficiency of the \KS can be improved by means of variable fixing \cite{guastaroba2014heuristic}.
Therefore, we shall illustrate a variant of the \KS based on a new strategy of variable fixing and apply it to solve large-scale \VMCPs. 
\begin{figure}[!htbp]
	\centering
	\includegraphics[width=\linewidth]{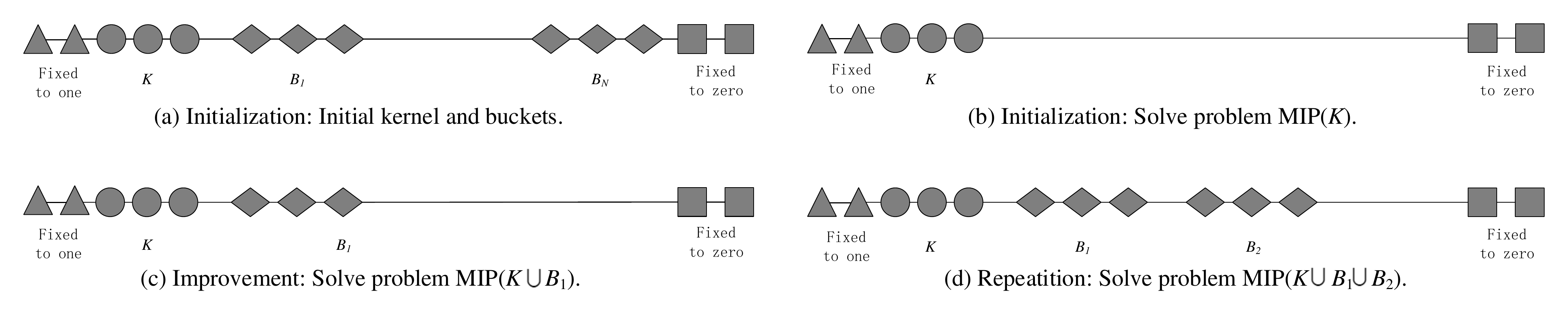}
	\caption{An illustrative example of the operation of the proposed \KS algorithm.}
	\label{fig:ksf_framework}
\end{figure}

Fig. \ref{fig:ksf_framework} describes an iteration in the proposed \KS algorithmic framework. 
The key change of our proposed \KS is 
to choose some binary and integer variables to fix their values, as compared to the standard \KS.
The strategies of fixed variables are crucial to the effectiveness and efficiency of our proposed \KS algorithm.
Therefore, we now first discuss the strategy of binary variable fixing, which reduces the number of binary variables in all restricted \MIP problems.
The straightforward strategies for fixing binary variables are as follows:
\begin{itemize}
	\item[1)] if $v_{j}^{*}=0$, then the associated binary variable is permanently fixed to zero; 
	\item[2)] if $v_{j}^{*}=1$, then the associated binary variable is permanently fixed to one;
\end{itemize}
where $v^{*}$, as stated, is the optimal solution for the linear relaxation of the original problem.
An existing variant of the \KS algorithm using strategies 1) and 2) was proposed in reference \cite{guastaroba2014heuristic}, 
and its numerical results indicated that the strategy improves the efficiency of the \KS algorithm with minor deteriorations of the solution quality. 
However, our preliminary experiments indicated that strategies 1) and 2) result in some \VMCP instances being infeasible.
For this reason, inspired by the basic linear programming theory \cite{Dantzig2003}, 
we use a more sophisticated but efficient strategy to avoid infeasibility occurring.
Each binary variable $v_{j}$ has an associated reduced cost $r_{j}^{*}$ value, which can be obtained by solving the \LP relaxation of the original problem.
The reduced cost is a lower bound on the increase of the \LP solution cost if the value of the variable is increased by one unit.
The strategies in our implementation are detailed as follows
\begin{itemize}
	\item[a)] if $r_{j}^{*} \geq \epsilon, ~\forall~j \in \mathcal{V}$ , then the associated binary variable $v_{j}$ is permanently fixed to zero 
	(see the squares in Fig. \ref{fig:ksf_framework}) and add it to set $\mathcal{Z}$;
	\item[b)] if $r_{j}^{*} \leq -\epsilon, ~\forall~j \in \mathcal{V}$, then the associated binary variable $v_{j}$ is permanently fixed to one 
	(see the triangles in Fig. \ref{fig:ksf_framework}) and add it to set $\mathcal{O}$.
\end{itemize}
where $\epsilon > 0$ that controls the number of fixed binary variables. 
In our implementation, we set $\epsilon=10^{-4}$.
Furthermore, since integer variables dominate the variables of the \VMCP, we apply the following similar strategy to fix integer variables, 
which further significantly improves the efficiency of our proposed \KS algorithm.
\begin{itemize}
	\item[c)] if $r_{j}^{*} \geq \epsilon, ~\forall~j \in \mathcal{G}$ , then the associated integer variable $g_{j}$ is permanently fixed to zero;
\end{itemize}

Our proposed \KS algorithm is summarized as in Algorithm \ref{alg:kernel_search_fix}. 
After permanently fixing the binary and integer variables according to strategies a), b), and c), 
and sorting the remainder of unfixed binary variables, 
we construct the initial kernel $\mathcal{K}$ and buckets sequence $\left\{\mathcal{B}_{i}\right\}_{i=1,\cdots,N}$.
Then we solve problem \MIP($\mathcal{U}$) where $\mathcal{U} \coloneqq \mathcal{K}$, 
but we find a few instances is infeasible due to the initial kernel is too small. 
To handle this situation, the proposed \KS follow \cite{guastaroba2017adaptive} to iteratively increase the size of the kernel until \MIP($\mathcal{U}$) is feasible.
A new kernel and bucket sequence are created when problem \MIP($\mathcal{U}$) becomes feasible.
The remainder steps are same to the standard \KS algorithm in Algorithm \ref{alg:kernel_search}. 

\begin{algorithm}[!htbp]
	\caption{The proposed \KS algorithm} 
	\label{alg:kernel_search_fix}
	\begin{algorithmic}[1]
		\Require \VMCP, $\bar{N}$, $T_{\rm max}$, $\epsilon$, and $\omega$.
		\Ensure $\text{UB}_{\min}$.
		\State Initialize $\text{UB}_{\min} \coloneqq +\infty$, $i \coloneqq  1$, $\mathcal{O}=\varnothing$, and $\mathcal{Z}=\varnothing$;
		\State Using the root node \LP relaxation of the original \VMCP instead of the \LP relaxation; 
		\State Using the following strategy to permanently fix binary/integer variable values:
		\State (a) if reduced cost $r_{j}^{*} \geq \epsilon, ~\forall~j \in \mathcal{V}$, the associated binary variable $v_{j}$ is permanently fixed to zero and $\mathcal{Z} \coloneqq \mathcal{Z} \cup \left\{v_{j}\right\}$;
		\State (b) if reduced cost $r_{j}^{*} \leq -\epsilon, ~\forall~j \in \mathcal{V}$, the associated binary variable $v_{j}$ is permanently fixed to one and $\mathcal{O}\coloneqq \mathcal{O} \cup \left\{v_{j}\right\}$;
		\State (c) if reduced cost $r_{j}^{*} \geq \epsilon, ~\forall~j \in \mathcal{G}$, the associated integer variable $g_{j}$ is permanently fixed to zero;
		\State Sort the unfixed binary variables in $\mathcal{V} \backslash \left(\mathcal{O} \cup \mathcal{Z}\right)$ according to a predefined sorting criterion; 
		\State Construct an initial kernel $\mathcal{K}$ by selecting the first $\vert \mathcal{K} \rvert$ binary variables in $\mathcal{V} \backslash \left(\mathcal{O} \cup \mathcal{Z}\right)$;
		\State Consider the binary variables belonging to $\mathcal{V} \backslash \left(\mathcal{K} \cup \mathcal{O} \cup \mathcal{Z}\right)$ and construct a sequence $\left\{\mathcal{B}_{i}\right\}_{i =1,\cdots,N}$ of buckets;
		\State Set maximum time limit $T_{i} \coloneqq T_{\rm max} / (N+1)$ of each bucket; 
		\State Solve problem \MIP($\mathcal{U}$) where $\mathcal{U} \coloneqq \mathcal{K}$;
		\If{no feasible solution to problem MIP($\mathcal{U}$) is found} \cite{guastaroba2017adaptive}
		\While{MIP($\mathcal{U}$) is not feasible} 
		\State Add the first $\vert \mathcal{K} \rvert \times \omega$  variables in the buckets sequence to set $\mathcal{U}$.
		\State Solve problem \MIP($\mathcal{U}$).
		\EndWhile
		\State Redefine the kernel $\mathcal{K} \coloneqq \mathcal{U}$ and buckets sequence $\left\{\mathcal{B}_{i}\right\}_{i=1,\cdots,N}$.
		\EndIf
		\While{$i \leq \bar{N}$}
		\State Add the two following constraints in \MIP($\mathcal{U} \cup \mathcal{B}_{i}$):
		\State (i) \eqref{eq:obj} $\leq {\rm UB}_{\rm min}$;
		\State (ii) $\sum\limits_{j \in B_{i}} v_{j} \geq 1$;
		\State Solve problem \MIP($\mathcal{U}$) where $\mathcal{U} \coloneqq \mathcal{U} \cup \mathcal{B}_{i}$ within time limit $T_{i}$ and denote its objective value by ${\rm UB}_{i}$; 
		\If{${\rm UB}_{i} < {\rm UB}_{\rm min}$}
		\State Update ${\rm UB}_{\rm min} \coloneqq {\rm UB}_{i}$;
		\EndIf
		\State Set $i \coloneqq i + 1$; 
		\EndWhile
	\end{algorithmic}
\end{algorithm}

\section{Numerical results}
\label{sect:numerical_results}
In this section, the effectiveness and efficiency of the proposed \KS algorithm for solving \VMCPs are evaluated by simulation experiments.
First, we perform computational experiments to demonstrate the performance of the proposed \KS algorithm 
for solving the \VMCPs, and compare it with the standard \MIP solver with three different settings. 
Then we evaluate the impact of the strategy of variable fixing on the performance of the \KS algorithm. 
All experiments were conducted on a cluster of Intel(R) Xeon(R) Gold 6140 @ 2.30GHz computers, 
with 192 GB RAM, running Linux (in 64 bit mode).

In our experiments, the proposed \KS was implemented in C++ linked with software 
IBM ILOG \cplex optimizer 20.1.0 \cite{CPLEX}.
After preliminary experiments, we set the following \cplex parameters in our proposed \KS.
To save computation time for obtaining root node \LP information, we decide not to apply the RINS heuristic (parameter RINSHeur) 
and \MIP heuristic (parameter HeurFreq), and turn off the feasiblility pump heuristic (parameter FPHeur) and local branching heuristic (parameter LBHeur). 
For all restricted problem \MIP($\mathcal{U}$), we choose the pseudo costs to drive the selection of the variable to branch on 
at a node (parameter VarSel), and generate mixed integer rounding cut (parameter MIRCuts) moderately.
All the other \cplex parameters were set to their default values. 
Following \cite{guastaroba2014heuristic}, the three settings for standard \MIP solver \cplex for comparison with the proposed \KS algorithm are as follows:
\begin{itemize}
	\item \CPXA: \cplex with all the parameters set to their default values with the exception that parameter MIPEmphasis that is set to feasibility.
	\item \CPXB: \cplex with all the parameters set to their default values with the exception that parameter RINSHeurStrategy is applied every 20 nodes.
	\item \CPXC: \cplex with all the parameters set to their default values with the exception that parameter LBHeurStrategy is turned on.
\end{itemize}
To gain insight the effectiveness of the proposed strategy of variable fixing, we consider three versions of the \KS algorithm: 
\begin{itemize}
	\item \KSFVG: our proposed \KS algorithm with variable fixing strategies a), b), and c). 
	\item \KSFV: has the same settings of \KSFVG with the exception that integer variables are not fixed (strategy c) is not applied).
	\item \KSF: standard kernel search algorithm without any strategy of variable fixing.
\end{itemize}
Finally, the time limit and the number of threads were set to 7200 (seconds) and 12. 
The optimal solutions for \VMCP instances are not available in the literature, so we use the best solution found by solver \cplex within 5 hours with the default setting, 
denoted as $f^{*}$, to validate the performance of three \cplex settings and three versions of the \KS algorithm.
\subsection{Testsets}
\label{sect:testsets}
All algorithms were tested on \VMCP instances with 5 \VM types and 10 server types with different features, 
as studied in \cite{mazumdar2017power}.
The \VM types and server types are reported in Tables \ref{table:first_dataset_vm} and \ref{table:first_dataset_server}, respectively. 

We generate four sizes of \VMCP instances, each having number $\lvert \mathcal{K} \vert \in \left\{7000,8000,9000.10000\right\}$ of servers. 
Each instance has the equal number of servers of each type.
In our test, the \VMCP instances are constructed as the following procedure.
First, we randomly select the element $k \in \mathcal{K}$ to subset $\hat{\mathcal{K}}$, with a probability $\alpha = 50\%$.
Then for each server $k \in \hat{\mathcal{K}}$, we iteratively assign a random number $\NIK$ of \VMs of type $i$ until
the maximum usage of the available resource load $\sigma_{k}$, defined by
\begin{equation}
\sigma_{k}=\max\biggl\{\frac{\sum_{i \in \ISET}\UIR \NIK}{\SKR}: \forall~r \in \RSET \biggr\}, \label{eq:generation_instance}
\end{equation}
exceeds a predefined value $\beta$.
Finally, to obtain parameter $\DINEW$ for all $i \in \mathcal{I}$, 
we iteratively assign a random number $\hat{d}^{\rm new}_{i}$ of \VMs of type $i$ until 
the maximum usage of the available resource load $\tau$, defined by
\begin{equation}
\tau=\max\biggl\{\frac{\sum_{i \in \ISET}\UIR \DINEW}{\sum_{j \in \KSET}\SKR}: \forall~r \in \RSET \biggr\}, 
\DINEW \coloneqq \DINEW + \hat{d}^{\rm new}_{i} \label{eq:generation_instance2}
\end{equation}
exceeds a predefined value $\gamma$.
In general, the larger $\beta$ and $\gamma$, the more \VMs and new incoming \VMs will be constructed.

As shown in references \cite{mazumdar2017power,gao2013multi,wu2016energy}, we consider the linear power consumption model as follows:
\begin{align} 
\label{eq:power_consumption_server}
P_{k}=P_{{\rm idle},k} +\left(P_{{\rm max},k}-P_{{\rm idle},k}\right)U_{k}.
\end{align}
$P_{{\rm idle},k}$ is the idle power consumption (at the idle state) of server $k$. 
$P_{{\rm max},k}$ is the maximum (peak) power consumption (at the peak state) of server $k$, and 
$U_{k}$ ($U_{k}\in \left[0,1\right]$) is the CPU utilization of server $k$.
Following \cite{mazumdar2017power}, the activation cost $\CRUN$ is set to $P_{{\rm idle},k}$,
the assignment cost $\CASSIGN$ is set to $\left(P_{{\rm max},k}-P_{{\rm idle},k}\right)\frac{u_{i,{\rm CPU}}}{s_{k,{\rm CPU}}}$,
and the idle power consumption $P_{{\rm idle},k}$ is set to 60\% of maximum power consumption $P_{{\rm max},k}$.
As illustrated in Section \ref{sect:problem_formualtion}, $\CMIG$ is set to $\CASSIGN$ in our experiments.

For each $ \lvert\mathcal{K}\vert \in \left\{7000,8000,9000,10000\right\}$, $\alpha=50\%$, $\beta \in \left\{20\%,40\%\right\}$, and $\gamma=50\%$,
50 \VMCP instances are randomly generated, leading to an overall 400 \VMCP instances.  

\begin{table}[!htbp]

	\centering
	\setlength{\tabcolsep}{10pt}
	\renewcommand{\arraystretch}{1.1}
	\caption{The five \VM types.}
	\label{table:first_dataset_vm}
	\begin{tabular}{cccc}
		\hline
		Type&CPU&RAM (GB)&Bandwidth (Mbps)\\
		\hline
		\VM 1&1&1&10 \\
		\VM 2&2&4&100\\
		\VM 3&4&8&300\\
		\VM 4&6&12&1000\\
		\VM 5&8&16 &1200\\
		\hline
	\end{tabular}

\end{table}

\begin{table}[!htbp]
	\centering
	\setlength{\tabcolsep}{1pt}
	\renewcommand{\arraystretch}{1.1}
	\caption{The ten server types.}
	\label{table:first_dataset_server}
	\begin{tabular}{ccccc}
		\hline
		Type&CPU&RAM (GB)&Bandwidth (Mbps)&Max power (W)\\
		\hline
		Server 1&4&8 &1000 &180\\
		Server 2&8&16 & 1000&200\\
		Server 3&10&16&2000&250\\
		Server 4&12&32&2000&250\\
		Server 5&14&32&2000&280\\
		Server 6&14&32&2000&300\\
		Server 7&16&32&4000&300\\
		Server 8&16&64&4000&350\\
		Server 9&18&64&4000&380\\
		Server 10&18&64&4000&410\\
		\hline
	\end{tabular}
	\vspace{-15pt}
\end{table}
\subsection{Efficiency of the proposed \KS algorithm}
\label{sect:efficiency_ksf}
In this subsection, we present numerical results to illustrate the efficiency of the proposed \KS algorithm 
compared with the standard \MIP solver with three different settings. 
\begin{table}[!htbp]
	\centering
	\setlength{\tabcolsep}{2pt}
	\renewcommand{\arraystretch}{1.2}
	\caption{Comparison results of \KSFVG and three \cplex settings}
	\label{table:comparison_results_1}
	\begin{tabular}{llllllllll}
		\hline
		$\left(\lvert \mathcal{K}\vert, \beta\right)$&\NVM & \multicolumn{2}{c}{\KSFVG}&\multicolumn{2}{c}{\CPXA}&\multicolumn{2}{c}{\CPXB}&\multicolumn{2}{c}{\CPXC} \\
		\cmidrule(r){3-4}\cmidrule(r){5-6}\cmidrule(r){7-8}\cmidrule(r){9-10}
		&&\GP&\TT&\GP&\TT&\GP&\TT&\GP&\TT \\
		\hline
		(7000,20\%)&24598& \textbf{0.0001}& \textbf{117.72}& 0.0003& 272.75  & 0.0002& 195.57  & 0.0002& 225.01 \\
		(7000,40\%)&27504 & \textbf{0.0001}& \textbf{116.82}& 0.0005& 226.81   & 0.0002& 197.91   & 0.0002& 219.55  \\
		(8000,20\%)&28095 & \textbf{0.0002}& \textbf{164.97} & 0.0006& 333.47   & 0.0002& 252.85   & 0.0002& 253.68  \\
		(8000,40\%)&31381 & \textbf{0.0001}& \textbf{161.02} & 0.0005& 298.13  & 0.0002& 225.37   & 0.0003& 285.85 \\
		(9000,20\%)&31549 & 0.0003& \textbf{199.82} & 0.0005& 461.56   & \textbf{0.0002}& 328.04   & 0.0004& 320.49 \\
		(9000,40\%)&35288  &\textbf{0.0002}& \textbf{174.24} & 0.0008& 392.13   & 0.0002& 277.90  & 0.0003& 349.02  \\
		(10000,20\%)&35077 &\textbf{0.0002}& \textbf{217.79} & 0.0007& 528.99  & 0.0003& 402.74  & 0.0003& 361.77 \\
		(10000,40\%)&39130 &\textbf{0.0001}& \textbf{197.75} & 0.0007& 449.78  & 0.0002& 346.64  & 0.0002& 408.90 \\
		\hline
		All & 31258&\textbf{0.0002} & \textbf{164.94} &0.0005 & 357.17 & 0.0002 & 269.79 & 0.0003 & 296.22  \\
		\hline
	\end{tabular}
\end{table}

In Table \ref{table:comparison_results_1}, we provide a summary of the computational results for the proposed \KS algorithm and  three \cplex settings.
For each data set, we report the average error (\GP) with respect to $f^{*}$ and the geometric mean of CPU time (\TT) in seconds.
The best solution values found by the proposed \KS algorithm or three \cplex settings are denoted as $f^{H}$.
The error for each instance is computed as $100\left(f^{H}-f^{*}\right)/f^{*}$, 
and then geometric averaged over all the instances belonging to the same data set to obtain statistic \GP.
As observed in Table \ref{table:comparison_results_1}, compared with \cplex settings \CPXA, \CPXB, and \CPXC, 
the CPU time of \KS algorithm taken by solving \VMCP is much smaller 
(169.94 seconds versus 357.17 seconds, 269.79 seconds, and 296.22 seconds).
In particular, we observe that the quality of the solutions found by the \KSFVG is not worse than \CPXA, \CPXB, and \CPXC, even slightly better.
From Table \ref{table:comparison_results_1}, we can conclude that the performance of \KSFVG 
is much better than three \cplex settings for all $\lvert\mathcal{K}\vert \in \left\{7000,8000,9000,10000\right\}$ and $\beta \in \left\{20\%,40\%\right\}$.
\begin{table}[!htbp]
	\centering
	\setlength{\tabcolsep}{5pt}
	\renewcommand{\arraystretch}{1.2}
	\caption{Worst error comparison results of \KSFVG and three \cplex settings}
	\label{table:comparison_results_2}
			\begin{tabular}{lllll}
				\hline
				$\left(\lvert\mathcal{K}\vert, \beta\right)$&\KSFVG&\CPXA&\CPXB&\CPXC \\
				&\WGP&\WGP&\WGP & \WGP\\
				\hline
				(7000,20\%)& \textbf{0.0071} & \textbf{0.0071}&0.0085&0.0082 \\
				(7000,40\%)& 0.0038&0.0065&\textbf{0.0034}&0.0095  \\
				(8000,20\%)& 0.0099&\textbf{0.0085}&\textbf{0.0085}&0.0088 \\
				(8000,40\%)& \textbf{0.0046}&0.0089&0.0069&0.0092 \\
				(9000,20\%)& \textbf{0.0078}&0.0084&0.0265&0.0100 \\
				(9000,40\%)&\textbf{0.0043}&0.0083&0.0573&0.0075 \\
				(10000,20\%)&\textbf{0.0063}&0.0163&0.0081&0.0067 \\
				(10000,40\%)&\textbf{0.0051}&0.0080&0.0070&0.0070\\
				\hline
				All  &\textbf{0.0099} &0.0163&0.0573&0.0100 \\
				\hline
			\end{tabular}
\end{table}

To gain more insight into the error of \KSFVG over three \cplex settings, we compare the worst error (\WGP) returned by \KSFVG and three \cplex settings.
Statistic \WGP shows the worst error calculated from all the instances that belong to the same data set.
The worst error comparison results of \KSFVG and three \cplex settings are summarized in Table \ref{table:comparison_results_2}.
From the table, we can conclude that the worst error returned by \KSFVG is less than 
that returned by three \cplex settings, especially for largest case ($\lvert\mathcal{K}\vert=10000$).

From the above computational results, we can conclude that the proposed \KSFVG algorithm is more efficient 
in solution quality and CPU time than the three \cplex settings for large-scale \VMCP instances. 

\subsection{Performance of the proposed strategy of variable fixing}
\label{sect:performance_comparison_kernel_search}
To address the advantage of applying the proposed strategy of variable fixing to the \KS algorithm, 
we compare \KSFVG with \KSFV and \KSF to solve \VMCPs.
The proposed \KSFVG has more fixed variables than \KSFV and \KSF does not have fixed variables.  

Table \ref{table:comparison_results_3} provides the computational results of the three versions of \KS algorithm. 
As expected, (i) the CPU time of \KSFVG is the least of the three versions of the \KS algorithm; 
(ii) the CPU time of \KSF is the most of the three versions of the \KS algorithm.
\begin{table}[!htbp]
	\centering
	\setlength{\tabcolsep}{4pt}
	\renewcommand{\arraystretch}{1.2}
	\caption{Comparison results of \KSFVG, \KSFV, and \KSF}
	\label{table:comparison_results_3}
	\begin{tabular}{lllllll}
		\hline
		$\left(\lvert\mathcal{K}\vert, \beta\right)$& \multicolumn{2}{c}{\KSFVG}&\multicolumn{2}{c}{\KSFV}&\multicolumn{2}{c}{\KSF}\\
		\cmidrule(r){2-3}\cmidrule(r){4-5}\cmidrule(r){6-7}
		&\GP&\TT&\GP&\TT&\GP&\TT\\
		\hline
		(7000,20\%)& 0.0001 & \textbf{117.72}& 0.0001&1143.62 & 0.0002& 2006.19 \\
		(7000,40\%)& 0.0001& \textbf{116.82}&0.0001& 781.32   & 0.0001& 1232.41   \\
		(8000,20\%)& 0.0002& \textbf{164.97} &0.0001& 1345.61  & 0.0001& 1964.78  \\
		(8000,40\%)& 0.0001& \textbf{161.02} &0.0001&  1111.64  & 0.0001& 1403.97  \\
		(9000,20\%)& 0.0003& \textbf{199.82} &0.0002& 1770.99  & 0.0002& 2538.79   \\
		(9000,40\%)&0.0002& \textbf{174.24} &0.0001& 1408.76  & 0.0002& 1732.06  \\
		(10000,20\%)&0.0002&\textbf{217.79}  &0.0003& 2480.27  & 0.0002& 3329.56 \\
		(10000,40\%)&0.0001& \textbf{197.75} & 0.0001& 1555.65& 0.0001& 2449.14  \\
		\hline
		All &0.0002 & \textbf{164.94}&0.0001 & 1376.21 & 0.0001 & 1988.61 \\
		\hline
	\end{tabular}
\end{table}
This is reasonable as the solution space for all restricted problem decreases with the number of fixed variables.
Furthermore, the quality of the solutions generally deteriorates with the number of fixed variables.
However, we observe that the quality of the solutions found by \KSFVG is slightly worse than \KSFV and \KSF, 
but the improvements in terms of CPU time are remarkable.
In some data sets, we observe that the average error of \KSFVG is even slightly less than that of \KSFV or \KSF.
This is due to the fact that some instances of \KSFV or \KSF cannot be solved within the time limit, resulting in their average error being larger.
\begin{table}[!htbp]
	\centering
	\setlength{\tabcolsep}{5pt}
	\renewcommand{\arraystretch}{1.2}
	\caption{Worst error comparison results of \KSFVG, \KSFV, and \KSF}
	\label{table:comparison_results_4}
	\begin{tabular}{llll}
		\hline
		$\left(\lvert\mathcal{K}\vert, \beta\right)$&\KSFVG&\KSFV&\KSF \\
		&\WGP&\WGP&\WGP\\
		\hline
		(7000,20\%)& 0.0071 &0.0061 & 0.0049\\
		(7000,40\%)& 0.0038&0.0025 &0.0016  \\
		(8000,20\%)& 0.0099&0.0051 &0.0061\\
		(8000,40\%)& 0.0046&0.0033 &0.0011\\
		(9000,20\%)& 0.0078&0.0048 &0.0078\\
		(9000,40\%)&0.0043&0.0070 &0.0048\\
		(10000,20\%)&0.0063&0.0079 &0.0059\\
		(10000,40\%)&0.0051&0.0046 &0.0057 \\
		\hline
		All &0.0099 &0.0079 &0.0078 \\
		\hline
	\end{tabular}
\end{table}
Next, we compare the worst error for \KSFVG, \KSFV, and \KSF, which is summarized in Table \ref{table:comparison_results_4}.
We observed that the worst error of the three versions of the \KS algorithm is greater with the decreasing value of $\beta$.
Similar behavior can be observed in the CPU time returned by Table \ref{table:comparison_results_3}.
In summary, we can conclude that the improvements of our proposed strategy of variable fixing in terms of CPU time for the \KS algorithm are significant as well as the impact on the quality of the solutions can be negligible.

\section{Conclusion}
\label{sect:conclusion}
\indent In this work, we have designed a new \KS algorithm for the solution of the large-scale \VMCP.
The proposed \KS algorithm is based on a new strategy of variable fixing, 
which is more efficient in terms of avoiding \VMCP infeasible than the existing \KS strategy, 
making it more suitable to solve large-scale \VMCPs.
Extensive computational experiments on large \VMCP instances show that our proposed \KS algorithm outperforms three different heuristic settings of the standard \MIP solver, and our proposed strategy of variable fixing significantly improves the efficiency of the \KS algorithm as well as the degradation of solution quality can be negligible.

\bmhead{Acknowledgments}
This work is supported by the National Natural Science Foundation of China (Grant no.: 12171052, 11971073, and 11871115). The computations were done on the high performance computers of State Key Laboratory of Scientific and Engineering Computing, Chinese Academy of Sciences.

\bibliography{temp}


\end{document}